\documentclass[12 pt]{article}
\usepackage{amsmath,amssymb,amsfonts}
\newtheorem{thm}{Theorem}[section]
\newtheorem{lem}[thm]{Lemma}
\newtheorem{prop}[thm]{Proposition}
\newtheorem{cor}[thm]{Corollary}
\newtheorem{rem}[thm]{Remark}
\newtheorem{ex}[thm]{Example}

\newtheorem{assu-def}[thm]{Assumption}

\newenvironment{pf}{\paragraph{Proof.}}{\par\medskip}
\newcommand{\qed}{\ifhmode\unskip\nobreak\fi\quad\ensuremath\diamond\par}

\newcommand{\pp}{\mathbb P}

\newcommand{\Pic}{\operatorname{Pic}}

\numberwithin{equation}{section}

\newcommand{\Si}{\Sigma}

\newcommand{\calC}{\mathcal{C}}

\newcommand{\calK}{\mathcal{K}}

\newcommand{\calM}{\mathcal{M}}

\newcommand{\calO}{\mathcal{O}}
\newcommand{\calP}{\mathcal{P}}

\newcommand{\bbP}{\mathbb{P}}

\begin{document}

\title{On surfaces with $p_g=q=2$ and non-birational bicanonical map \thanks{2000 Mathematics
Subject Classification: 14J29}}
\author{Ciro Ciliberto, Margarida Mendes Lopes}
\date{}
\maketitle
\begin{abstract}  The present
paper is devoted to the classification of irregular surfaces
of general type with $p_g=q=2$ and non birational
bicanonical map. The main result is that, if $S$ is such a
surface and if $S$ is minimal with no pencil of curves of
genus $2$, then   $S$ is a double cover of a
principally polarized abelian surface $(A,\Theta)$, with
$\Theta$ irreducible. The double cover $S\to A$ is branched
along a divisor $B\in |2\Theta|$, having at most
double points and so $K_S^2=4$.

\end{abstract}
\section{Introduction}

If a smooth surface $S$ of general type has a pencil of curves of genus $2$, i.e. it has a
morphism to a curve whose general fibre $F$ is a smooth irreducible curve of genus 2, then the line
bundle ${\calO}_S(K_S)\otimes {\cal O}_F$ is the canonical bundle on $F$, and therefore the bicanonical
map $\phi$ of $S$ cannot be birational. Since this property is, of course, of a birational nature, the same remark applies
if $S$ has a rational map to a curve whose general fibre is an irreducible curve with geometric genus $2$.\par

We call this exception to the birationality of the bicanonical map $\phi$
 the {\it standard case}. A {\it non-standard case} will be the one of a surface of
general type $S$ for which $\phi$ is not birational, but there is no pencil of curves of
genus 2.
The classification of the non-standard cases has a long history and we refer  to the
expository paper \cite{ci} for information on this problem. We will just mention here the fact that the non-standard
cases with
$p_g\geq 4$ are all regular.

\par The classification of non-standard irregular surfaces has been considered
by Xiao Gang in \cite{xiaocan} and by F. Catanese and the authors of the present paper in \cite{ccm}. Xiao Gang
studied the general problem of classifying the non-standard cases  by taking the point of view of the projective
study of the image of the bicanonical map. The outcome of his analysis is a list of numerical
possibilities for the invariants of the cases which might occur. More precise results have been obtained in
\cite{ccm}, where the first significant case $p_g=3$ has been considered. Indeed in \cite {ccm} it is shown, among
other things, that  a minimal irregular surface $S$ with
$p_g=3$ presents the non-standard case  if and only if  $S$
is isomorphic to the symmetric product of a smooth irreducible curve of
genus $3$, thus $p_g=q=3$ and $K^2=6$.\par

 In the present paper we study this problem for surfaces with
$p_g=q=2$  and we prove the following result,  which rules out a substantial number of possibilities
presented in \cite{xiaocan}:

\begin{thm}\label{teorema}  Let $S$ be a minimal surface of general type with
$p_g=q=2$. Then $S$  presents the non-standard case if and only if $S$ is a double cover of a
principally
polarized abelian surface $(A,\Theta)$, with $\Theta$ irreducible. The double cover $S\to A$ is branched along
a symmetric divisor $B\in |2\Theta|$, having at most double points. One has $K_S^2=4$.\end{thm}

Surfaces with $p_g=q=2$ are still far from being understood. The list of known examples of surfaces of general type with
$p_g=q=2$ is relatively small (see \cite{zucconi1}, \cite{zucconi2}) and there are several constraints for their
existence known. Here we  only mention that there are various restrictions for the existence of a genus $2$ fibration (see
\cite{xiao}) and also that M. Manetti,  working on the Severi conjecture, showed in particular that if $p_g=q=2$, $K_S$ is
ample and $K_S^2=4$ then $S$ is a double cover of its Albanese image (see \cite{man}).

To prove our classification theorem \ref {teorema} we first show that the degree of the bicanonical map
is $2$ for surfaces presenting the non-standard case, then we study the possibilities for the quotient surface by the
involution induced by the bicanonical map, and finally we show that the unique case which really occurs is the one described
above. We use a diversity of techniques, which may be useful in other contexts.
\par

The paper is organized as follows. In section \ref{properties} we list  the properties of surfaces $S$ with $p_g=q=2$
that we need. In section \ref{K9} we characterize, by a small adaptation of a proof in \cite{ccm}, the surfaces
$S$  presenting the non-standard case with
$K_S^2=9$, and in particular we verify that there is no such surface with $p_g=q=2$. In section
\ref{paracanonical} we establish some properties of the paracanonical system and then we use these results in
section \ref{deg} to conclude that for the non-standard cases $S$ with $p_g=q=2$ the degree of the bicanonical map
is $2$. Thus there is an involution
$i$ induced by the bicanonical map on $S$. We consider the quotient surface $\tilde\Sigma:=S/<i>$   and the
projection map $p:S\to
\tilde\Sigma$. In section
\ref{sete}  we discuss the various possibilities for $\tilde\Sigma$,
 showing that the only one which can really occur is that $\tilde \Sigma$ is a minimal surface of general type  with
$p_g(\tilde\Sigma)=2$,
$q(\tilde\Sigma)=0$,
$K_{\tilde\Sigma}^2=2$ and with 20 nodes. Moreover we show  that the double cover $p$ ramifies exactly over the $20$
nodes. Finally in section \ref{main}, using this description, and some results on Prym varieties contained in
\cite{cpt}, we finally prove theorem \ref{teorema}.
\medskip

\paragraph{Acknowledgements} The present collaboration takes place in the framework of the european contract EAGER, no.
HPRN-CT-2000-00099. The second author is a member of CMAF
and of the Departamento de Matem\'atica da Faculdade de Ci\^encias da Universidade
de Lisboa.

We are indebted to Rita Pardini for interesting discussions on the subject of this paper and in
particular for having pointed out the use of the Kawamata--Viehweg vanishing theorem in section \ref{sete}.

We thank prof. Fabrizio Catanese for having communicated to us some of his ideas on this
subject.

We  dedicate this paper to the memory of Paolo Francia, with whom we first started  on
this subject.

\paragraph{\bf Notations and conventions} We work over the complex numbers. All
varieties are assumed to be compact and algebraic. We do not distinguish
between line bundles and divisors on a smooth variety, using the additive
and the multiplicative notation interchangeably. Linear equivalence is denoted
by $\equiv$ and
numerical equivalence by $\sim$. A {it node} on a surface is an ordinary double point (i.e a singularity of type
$A_1$). The exceptional divisor of a minimal desingularization of a node is a rational irreducible curve $A$ with $A^2=-2$,
usually called a
$(-2)$-curve.

As already mentioned, we will say that a surface $S$ of general type {\it presents the non-standard case}, or that it is a
{\it non-standard case}, if $S$ has no pencil of curves of geometric genus $2$ and the bicanonical map of $S$ is not
birational.

 The remaining notation is standard in
algebraic geometry.

\section{ Some properties of surfaces
with $p_g=q=2$.}\label{properties}

The minimal
surfaces $S$ of general type  with $p_g=q=2$ have various interesting properties (cf. \cite{zucconi1}, \cite{zucconi2}). In
this section we only mention those that we will need further on.

\begin{prop}\label{p2}  Let $S$ be a minimal
surface of general type with $p_g=q=2$. Then:

\vskip0.2truecm\noindent i) $4\leq K_S^2\leq 9$;

\vskip0.2truecm\noindent ii) if $K_S^2=8,9$, there are no rational smooth curves on $S$ (in particular
${\calO}_S(K_S)$ is ample), and, if
$K_S^2= 7$ and
${\calO}_S(K_S)$ is not ample, then
$S$ contains either one irreducible
$(-2)-$curve, or two forming an $A_2$ configuration. Furthermore if $K_S^2=9$, $S$ does not contain  elliptic
 curves.

\end{prop}

\begin{pf}
 i)
The
first inequality follows from the  inequality $K^2\geq 2p_g$ for minimal
irregular surfaces (see \cite{De}), and the second from the
inequality
$K^2\leq 3c_2$.\par

ii) follows from Miyaoka's and Sakai's inequalities (see \cite{miyaoka} and \cite{sa}) for the number of
rational or elliptic curves on a non ruled minimal surface.
\qed
\end{pf}

\medskip
We will need also to consider the Albanese image of these surfaces. First we recall the following facts which we will use
repeatedly:
\begin{lem}\label{fibration}
{\rm (see  \cite{Be}, pg. 343; \cite{bpv}, pg. 97) } Let $S$ be a minimal surface and let $f:S\to B$ be a genus $b:=g(B)$ pencil of
curves of genus $g\geq 2$. Then

  \vskip0.2truecm\noindent i) $K_S^2\geq 8(g-1)(b-1)$;\par
\vskip0.2truecm\noindent ii) $c_2(S)\geq 4(g-1)(b-1)$
 \vskip0.2truecm\noindent iii) $q\leq g+b$.

\par Furthermore if equality holds in i) then the curves of the pencil have
constant modulus, if  equality holds in ii) every fibre of $f$ is smooth and if  equality holds in iii)  $S$ is birationally
equivalent to a product of $B$ with the general fibre of $f$.\end {lem}

  Using this lemma we obtain the following:

\begin{prop}\label{alb} Let $S$ be a minimal
surface of general type with $p_g=q=2$ for which the  Albanese morphism $a: S\to
A:=Alb(S)$ is not surjective. Then $a(S)=B$ is a genus $2$ curve, the Albanese pencil $a: S\to B$
has smooth, connected fibres $F$ of genus $2$ with constant modulus  and $K_S^2=8$.

\end{prop}
\begin{pf} Since $q(S)=2$, the Albanese image of $S$ is a genus $2$ curve $B$. Then
the remainder of the assertion is a consequence of    lemma
\ref{fibration} and $\chi({\calO}_S)=1$.
\qed  \end{pf}

\begin{cor}\label{defranchis}  Let $S$ be a minimal
surface of general type with $p_g=q=2$. If  $\omega, \omega'$ are
two 1-forms which generate
$H^0(S,\Omega^1_S)$ and $\omega\wedge\omega'\equiv 0$, then the  Albanese morphism $a: S\to
A:=Alb(S)$ is not surjective, the Albanese pencil $a: S\to B$
has smooth, connected fibres $F$ of genus $2$ with constant modulus  and $K_S^2=8$.
\end{cor}
\begin{pf}  The assertion follows from the theorem of  Castelnuovo-De Franchis  (see
e.g.
\cite{bpv}, pg. 123) and the previous proposition.
\qed  \end{pf}

\vskip0.2truecm\noindent Finally we notice that, if the  surface $S$ of general type  with $p_g=q=2$ has
a genus $2$ fibration,  then  the canonical system is not
composed with the genus 2 fibration (see \cite{xiao}, Theorems 2.1, pg. 16, and Theorem 5.1, pg. 71). As a consequence we
have:

\begin{prop}\label{irreducible} Let $S$ be a minimal
surface of general type with $p_g=q=2$ and write $|K_S|=|M|+Z$, where $|M|$ is the moving part of $|K_S|$ and $Z$ the
fixed part. Then the general curve in $|M|$ is irreducible.
\end{prop}

\begin{pf} Assume otherwise. Then $|M|$ is composed with an irrational pencil $\calP$. If $F$ is a generic fibre of
$\calP$, $|M|=aF$ where $a\geq 2$, and furthermore $F^2=0$. Since $F$ is not a genus 2 curve,
$K_S\cdot F\geq 4$. Since $K_S^2\leq 9$, we see that either $K_S\cdot Z=1, K_S^2=9$ or $K_S\cdot Z=0,
K_S^2=8$. This cannot occur. Indeed, in the former case
$S$ would contain a curve $\theta$ with $K_S\cdot \theta=1$, hence $\theta$ would be rational or elliptic, whereas in the latter
case $S$ would contain a $(-2)-$curve. In either case we would have a  contradiction to proposition \ref{p2}, ii).
\qed  \end{pf}

\section{The case $K_S^2=9$.}\label{K9}

In \cite{reider} I. Reider proved that if $S$ is a minimal surface of general type
with $K_S^2\geq 10$ and the bicanonical map is not birational, then $S$
presents the standard case. In Proposition (1.1) of \cite{ccm},
it is proven that the same holds if $K_S^2=9$ and $p_g\geq 3$, unless $p_g=6, K_S^2=9$,
and $S$ is the Du Val-Bombieri's surface described in \cite{duval} and in \cite{bombieri}, pg. 193. In fact this
result can be extended:

\begin{prop} \label{9} Let $S$ be a minimal
surface of general type with $K_S^2=9$ such that the bicanonical
map is not birational. Assume that $S$ presents the non-standard case. Then $p_g=6$,
$q=0$ and $S$ is the  Du Val-Bombieri's surface.
\end{prop}

\begin{pf}  To prove the assertion it suffices to use the proof  of
proposition (1.1) of \cite{ccm}.  There, the
assumption $p_g\geq 3$ is only necessary for the proof of Claim $4$. But  Claim $4$ can be proved without using the
assumption  on $p_g$. In fact, since $K_S-D\sim 2D$ is big and nef, Mumford's vanishing
theorem (see \cite{mumford}, pg. 250), yields  $h^1(S,{\calO}_S(2K_S-D))=0$. Thus the   map
$H^0(S,{\calO}_S(2K_S))\to H^0(D,{\calO}_D(2K_S))$ is surjective, which in turn implies that $D$ is hyperelliptic.
\qed  \end{pf}

\section{ The paracanonical system in the case $p_g=q=2$.}\label{paracanonical}

Let $S$ be a minimal irregular surface of general type.  If $\eta\in Pic^0(S)$ is a point, we can
consider the linear system $|K_S+\eta|$. A curve in $|K_S+\eta|$ is a {\it paracanonical curve} on $S$.
\par

Assume that the  Albanese image of $S$ is a surface. Given a general point
$\eta\in Pic^0(S)$, one has,  by \cite{gl}, Thm. 1.,  $h^1(S,\calO_S(\eta))=0$ and $dim
|K_S+\eta|= \chi ({\calO}_S) -1$.

For $\eta\in Pic^0(S)$, let $C_\eta$ be the general
curve in $|K_S+\eta|$.  The curves $C_\eta$
describe, for
$\eta\in Pic^0(S)$ a general point, a continuous system ${\calK}$ of curves on
$S$, of dimension
$q+dim |K_S+\eta|= q+\chi ({\calO}_S) -1=p_g$. This is what we will call the {\it
main paracanonical system} of $S$. \par

Assume now that $S$ is a minimal surface of general type with $p_g=q=2$, for which the Albanese map
$a: S\to A:=Alb(S)$ is surjective. The main paracanonical system of $S$ has dimension $2$
and, for
$\eta\in Pic^0(S)$ a general point, the curve $C_\eta\in |K_S+\eta|$ is linearly
isolated. We  write $C_\eta=F+M_\eta$, where $F$ is the fixed part of the
continuous system
${\calK}$ and $M_\eta$ the movable part, and we denote by ${\calM}$ the continuous,
$2$-dimensional system described by the curve $M:=M_\eta$. This system is parametrized by a surface $P$ which is birational to
$Pic^0(S)$.
\par
\begin{lem}\label{paracanonicalcurve} Let $S$ be a minimal
surface of general type with $p_g=q=2$ presenting the non-standard case. Let $C_\eta=F+M_\eta$ be
the general paracanonical curve. Then either:\par

\noindent (i) $M:=M_\eta$ is irreducible and $M^2\geq 3$, or;\par

\noindent (ii) $F=0$ and $M$ is reducible as $M=M_1+M_2$, with $M_1$ and $M_2$
irreducible each varying in two $1$-dimensional systems of curves ${\calM}_1$, ${\calM}_2$. The following possibilities can
occur:

$$ (a)\quad M_1^2=M_2^2=0, M_1\cdot M_2=4, K_S^2=8$$

$$(b)\quad M_1^2=M_2^2=M_1\cdot M_2=2, M_1\sim M_2, K_S^2=8.$$

\end{lem}

\begin{pf} Suppose that $M$ is irreducible. Then $M^2>0$, otherwise
${\calM}$ is a pencil, whereas we know it has dimension $2$. The case $M^2=1$ is excluded by
proposition (0.14, iii) of \cite{ccm}. The case $M^2=2$ is also excluded by theorem (0.20) of \cite{ccm}. This
proves (i).\par

 Suppose that $M$ is reducible. Since ${\calM}$ is a two dimensional system parametrized by $Pic^0(S)$,
  $M$ must consist of two distinct irreducible components $M=M_1+M_2$.\par

Suppose $M_i^2=0$ for one of $i=1,2$. Then $M_i$ varies in a pencil ${\calM}_i$ of curves of genus at
least $3$ and so $K_S\cdot M_i\geq 4$. If instead $M_i^2>0$, then, by proposition (0.18) of \cite{ccm},
$M_i^2\geq 2$ and one has again $K_S\cdot M_i\geq 4$, by the 2-connectedness of the paracanonical curves. In both
cases
$$K_S^2=K_S\cdot F+K_S\cdot M_1+K_S\cdot M_2\geq 8,$$ and, so, by proposition \ref{9}, one has $K_S^2=8$ and
$K_S\cdot M_1=K_S\cdot M_2=4$,
$K_S\cdot F=0$. Since $S$ does not contain $(-2)-$curves, one has $F=0$ and we have the two numerical possibilities listed in
(ii). \qed
  \end{pf}
\bigskip

\begin{lem}\label{con}  Let $S$ be a minimal
surface of general type with $p_g=q=2$ presenting the non-standard case and let $C_\eta=F+M_\eta$ be as in case i) of
lemma
\ref  {paracanonicalcurve}.  Then:

\par\noindent i) if $F\neq 0$, then   $F\cdot M_\eta=2$ and $F$ is
$1$-connected. \par
\par\noindent ii) if $F\neq 0$ and $\eta$ is  general,  the image of the restriction map
$$H^0(S,{\calO}_S(2K_S))\to H^0(M_\eta,{\calO}_{M_\eta}(2K_S))$$ has at most codimension $1$ in
$H^0(M_\eta,{\calO}_{M_\eta}(2K_S))$.
\end{lem}

\begin{pf} i) Let $M:=M_\eta$. If $F\cdot M=2$, the 2-connectedness of the canonical divisors and lemma (A.4)
of \cite{cfm} imply that $F$ is 1-connected.
 To show that $F\cdot M=2$ first we claim that $F\cdot M\leq 4$. Indeed, proposition
\ref{9} yields  $K_S^2\leq 8$  and  lemma \ref{paracanonicalcurve}, i) yields  $M^2\geq 3$. Therefore

  $$8\geq K_S^2\geq
K_S\cdot M=M^2+F\cdot M\geq 3+F\cdot M.$$

So $F\cdot M$ being even implies  $F\cdot M\leq 4$.  \par
Now we show that $F\cdot M= 4$ cannot occur. Suppose otherwise. Then  from $8\geq K_S^2=M^2+8+F^2$
and
$K_S\cdot M=M^2+M\cdot F\geq 7$ we have the possibilities:\medskip

\par\noindent a) $K_S^2=7$, $K_S\cdot F=0$, $F^2=-4$,
\par\noindent b) $K_S^2=8$,
 $K_S\cdot F=1$, $F^2=-3$
\par \noindent c) $K_S^2=8$, $K_S\cdot F=0$, $F^2=-4$.\medskip

 The first possibility implies that $F$ contains two
disjoint  $(-2)$-curves, whilst the second and third imply that $F$ contains a smooth rational curve. This is impossible
by proposition \ref{p2}, ii). Therefore $F\cdot M=2$ and so $F$ is 1-connected.

\medskip

\noindent ii) Note that, since $H^1(S,{\calO}_S(2K))=0$, the codimension of  the image of the restriction map
$$H^0(S,{\calO}_S(2K))\to H^0(M_\eta,{\calO}_{M_\eta}(2K))$$ is exactly $h^1(S,{\calO}_S(K_S+F-\eta))$, which by duality
is equal to $h^1(S,{\calO}_S(\eta-F))$.
Consider the exact sequence $$0\to {\calO}_S(\eta-F)\to {\calO}_S(\eta)\to {\calO}_F(\eta)\to 0$$ which yields the
long exact sequence $$0\to H^0(S,{\calO}_S(\eta-F))\to H^0(S,{\calO}_S(\eta))\to H^0(F,{\calO}_F(\eta))\to$$
$$\to H^1(S,{\calO}_S(\eta-F))\to H^1(S,{\calO}_S(\eta))\to...$$
Since $h^0(S,{\calO}_S(\eta))=0$ and  $h^1(S,{\calO}_S(\eta))=0$, for $\eta$ general (by \cite{gl}, Thm. 1), we see that
$h^1(S,{\calO}_S(\eta-F))=h^0(F,{\calO}_F(\eta))$. Now ${\calO}_F(\eta)$ has degree $0$ on every component of
$F$. By  the first part  of the lemma, $F$ is 1-connected and so by  corollary (A.2) of \cite{cfm},
$h^0(F,{\calO}_F(\eta))\leq 1$ (with equality holding if and only if $ {\calO}_F(\eta)\simeq {\calO}_F$).
\qed  \end{pf}

\section{ The degree of the bicanonical map}\label{deg}

In the present section we  prove the following result:
\begin{prop}\label{main5} Let $S$ be a minimal
surface of general type with $p_g=q=2$. Assume
that $S$ presents the non-standard case. Then the degree $\sigma$ of the bicanonical map is $2$.

\end{prop}
\begin{rem}  {\rm For completeness let us point out that if $S$ has a genus $2$ fibration then  the degree $\sigma$ of the
bicanonical map is either $2$ or $4$, ( see
\cite{xiao}) and $\sigma=4$ does occur (cf. remark \ref{deg4}).}
\end{rem}
First of all we treat the case $K_S^2=8$, adapting a proof which appears in \cite{mp}.

\begin{prop}\label{k8} Let $S$ be a minimal
surface of general type with $p_g=q=2$ and $K_S^2=8$  presenting the non-standard case. Then the degree $\sigma$ of
the bicanonical map is
$2$.
\end{prop}
\begin{pf}  Let $\phi$ be the bicanonical map of $S$. Notice that $(2K_S)^2=4K_S^2=32$
and
$h^0((S,{\calO}_S(2K_S))=K_S^2+1=9$. Then the degree of $\Sigma=:\phi(S)$ is $\frac{32}{\sigma}\geq 7$, hence
$\sigma$ is either $2$ or
$4$.\par

Suppose $\sigma=4$.  In
this case $\Sigma$  is a surface of degree $8$ in $\bbP^8$. The list of such surfaces is known (see  \cite {na},
Thm. 8). Since
$|2K_S|$ is a complete linear system, $\Sigma$ can be one of the following:\medskip

\par\noindent a) the Veronese embedding in $\bbP^8$ of a quadric  in
$\bbP^3$;

\par\noindent b) a Del Pezzo surface, i.e. the image of $\bbP^2$ by the
rational map associated to the linear system
$|3L\otimes {\cal I}_{x|\bbP^2}|$, where $L$ is a line and
$x$ is a point of $\bbP^2$;

 \par\noindent c) a cone over an elliptic curve of degree $8$ in $\bbP^7$.
\vskip0.2truecm\par

We are going to prove the result by showing that none of these cases can occur. \par First we
consider case c).  Take the pull back $F$ of a line in the cone. Then $2K_S\cdot F=4$, hence
$K_S\cdot F=2$. The  index theorem then yields
$F^2=0$, and therefore we would have a genus $2$ pencil on $S$.\par

 In case a) $2K_S\equiv 2H$, where $H$ is the pull back of the hyperplane section of
$\Sigma$. Then $\eta=H-K_S$ is a nontrivial 2-torsion element in $\Pic S$, since
$p_g(S)=2$ whereas $h^0(S, {\calO}_S ( K_S+\eta)=4$. The \'etale double cover $\pi\colon Y\to S$ given by
$2\eta\equiv 0$ has invariants $\chi (Y)=2$, $K_Y^2=16$.
In addition
$p_g(Y)=p_g(S)+h^0(S,{\calO}_S(K_S+\eta))=6$   so that
 $q(Y)=5$. Then, since $q(S)=2$, the subspace $V^-$ of $H^0(Y,\Omega^1_Y)$ containing the
antiinvariant
$1-$forms by the involution $\iota$ determined by $\pi\colon Y\to S$  has dimension $3$. Since the image of $\wedge^2 V^-$ in
$H^0(Y,\Omega^2_Y)$ is contained in the subspace of invariant 2-forms which is 2-dimensional, we conclude  that there are two
independent $1-$forms $\omega,\omega'\in V^-$ such that
$\omega\wedge \omega'\equiv 0$ and so by the theorem of Castelnuovo-De Franchis there exists a fibration $g:Y\to B$ with
$b:=g(B)\geq 2$ (cf. also
\cite{cetraro}, corollary (4.8)).

 Let
$f$ be the genus of a general fibre $F$ of
$g$. Suppose $f=2, b\geq 3$. Then the curve $F'=\iota(F)$ cannot dominate $B$ via $g$. Hence $F'$ is again a curve of the
pencil  $g:Y\to B$. It cannot be the case that $F'=F$, otherwise $\pi(F)$ would be a moving curve of genus $0$ or $1$ on $S$, a
contradiction. In conclusion $F\neq F'$ and $\pi(F)=\pi(F')$ would be a curve of genus 2 on $S$ varying in a pencil, a
contradiction.

Now, by lemma \ref{fibration}, i),   we have $K_Y^2=16\geq 8 (f-1)(b-1)$ and, by lemma
\ref{fibration}, ii), $5=q(Y)\leq f+b$. This forces $f=3$, $b=2$ or viceversa and so  $Y$ is birational to $B\times F$
(see again lemma
\ref{fibration}, ii)). Hence
$Y$ has a pencil of curves of genus $2$, whose image on $S$, by what we observed above, is again a genus $2$ pencil, against
our hypothesis. Thus also case a) does not occur.\par

Finally we consider case b).  We abuse notation and we denote by $L$  the image on $\Sigma$ of a line of $\bbP^2$. Let
$2L+L_0$ be the hyperplane section of $\Sigma$.  We have
$2K_S\equiv \phi^*(2L+L_0)$, and so $\phi^*(L_0)\equiv 2(K_S-\phi^*L)$. \par

Choose
$L_0$  such that $\phi^*(L_0)$ is a smooth irreducible curve  and consider the double cover $Y$ of $S$ branched over
$\phi^*(L_0)$ and determined by
$K_S-\phi^*(L)$. The double cover formulas give $\chi(Y)=3$, $K_{Y}^2=24$, $p_g(Y)=p_g(S)+h^0(S, {\cal
O}_S(2K_S-\phi^*L))=7, $ so that $q(Y)=5$.

 Notice that $|\phi^*(L_0)|$ is a genus 3
pencil on
$S$. The pull back of it to $Y$ is either a  rational pencil of curves of genus 5,   or a genus 3 pencil. In the former case
$Y$ would be birational to the product of $\bbP^1$ by a curve of genus $5$ (see again lemma
\ref{fibration}), which is not
possible. In the other case let $b$ be the genus  of the base curve of the pencil. As before $b\geq 2$, because
$b+3\geq q(Y)=5$. On the other hand lemma
\ref{fibration}, i)  yields $K_{Y}^2=24\geq 16(b-1)$. Hence $b=2$ and as above we conclude
that $Y$ is birational to  a product of a genus 2 and a genus 3 curve, which is impossible because
$p_g(Y)=7$.
\qed
\end{pf}
\bigskip

Before continuing towards the  proof of proposition \ref{main5}  we need to recall some facts about continuous systems of
curves on a surface. For the basic definitions, we defer the reader to \cite {ccm}, \S 0.
Given an irreducible, continuous system ${\calC}$ of curves of dimension $r$ on a surface $S$, the {\it
index} $\nu:=\nu_{\calC}$ of ${\calC}$ is the number of curves of ${\calC}$ passing
through $r$  general points of $S$. Of course  $\nu\geq 1$.  A system ${\calC}$ is called an {\it involution} if its index
is $\nu=1$. Typical examples of involutions are:\par
\noindent (i) the linear systems;\par
\noindent (ii) pencils, or, more generally systems {\it composed with pencils}. This means that there
is a pencil $f: S\to B$ and an involution of divisors on $B$ such that the curves of ${\calC}$ are
pull-back, via $f$ of divisors of an involution on $B$.\par

The classical theorem of Castelnuovo-Humbert tells us that these are essentially the only involutions.

\begin{thm}[Castelnuovo-Humbert] {\rm (see \cite{chiantini}, \S 5)}  \label{Humbert}   Let $S$ be
a smooth, irreducible, projective surface and let ${\calC}$ be an
$r$-dimensional involution on $S$ which has no fixed divisor and whose general divisor $C$ is reduced.
Then either ${\calC}$ is a linear system or it is composed with
a pencil.
\end {thm}

We will use this theorem to prove the following basic result:

\begin{prop} \label{grado} Let $S$ be a minimal
surface of general type with $p_g=q=2$. Assume
that $S$ presents the non-standard case. Let $C_\eta=F+M_\eta$ be the general paracanonical curve and suppose
that $M:=M_\eta$ is irreducible. Then the restriction of the bicanonical map $\phi$ to $M$ is a
birational map of $M$ onto its image.\end{prop}

\begin{pf} First we consider the case  $F=0$. Then the
arithmetic genus $g$ of $M$ is $g=K_S^2+1$. Since, by \cite{gl}, thm. 1, $h^1(S,{\cal O}_S(K_S+\eta))=0$ for
 a general point $\eta\in Pic^0(S)$,  $|2K_S|$ cuts out on $M=M_\eta$ a non-special,
base point free complete $g^{g-2}_{2g-2}$. We will argue by contradiction and we will suppose from now on
that this series is composed with an involution $\tau:=\tau_M$ of degree $\delta\geq 2$. Then we must
have $2g-2\geq \delta(g-2)$ which yields $\delta\leq 2+ \frac{2} {g-2}=2+\frac {2} {K_S^2-1}$. Since, by
proposition \ref{p2}, one has $K_S^2\geq 4$, we see that $\delta=2$. This means that $\phi(M)$ is a
linearly normal curve of degree $g-1$ in $\bbP^{g-2}$, whose arithmetic genus is $1$. Notice that two
distinct points $x,x'$ are conjugated in $\tau_M$ if and only if $\phi(x)=\phi(x')$. \par

\vskip0.3truecm\noindent  {\it Claim 1: Let $M, M'$ be general curves in ${\cal M}$, then $M\cap M'$
does not contain four distinct points $x,y,x', y'$ such that $\phi(x)=\phi(x')$ and $\phi(y)=\phi(y')$.
}\vskip0.3truecm

Otherwise we would have $h^0(M, {\cal O}_M(M'))\geq h^0(M, {\cal
O}_M(x+x'+y+y'))=2$. On the other hand, since $h^1(S,{\calO}_S(\eta))=0$ for
$\eta\in Pic^0(S)$ a general point, $|M'|$ cuts out a complete linear series on $M$. Since $M'$ is
linearly isolated, we find a contradiction. \vskip0.3truecm

Let $x$ be a point on $S$. We denote by ${\cal M}_x$ the
system of curves in ${\cal M}$ passing through $x$.

\vskip0.3truecm\noindent  {\it Claim 2: Let $x$ and $x'$ be general points on $M$
conjugated in $\tau$, i.e. such that $\phi(x)=\phi(x')$. Every irreducible component of ${\cal M}_x$ is
a $1$-dimensional system of curves. Consider the union  of all of these components containing $M$. Every
curve in such a union contains $x'$. } \vskip0.3truecm

Let $M"$ be the general curve in a component ${\cal M}'$ of the union in question and let $x_{M"}$ be the
point conjugated to $x$ in the involution $\tau_{M"}$ on $M"$. Since $\phi(x_{M"})=\phi(x)$ and $\phi$ is
generically finite, then $x_{M"}$ belongs to a finite set when $M"$ varies in ${\cal M}'$, and therefore it
stays fixed when $M"$ varies in ${\cal M}'$. Since $x_M=x'$ we have $x_{M"}=x_M=x'$, proving the claim.
\vskip0.3truecm

It is appropriate to denote by ${\cal M}_{M,x,x'}$ the union of all components of ${\cal M}_x$
containing $M$. Since ${\cal M}$ is parametrized by a surface $P$ birational to $Pic^0(S)$, the system
${\cal M}_{M,x,x'}$ corresponds to a reduced curve $D_{M,x,x'}$ on $P$. This curve might be
reducible, but all of its irreducible components pass, by definition, through the point $m$ of $Pic^0(S)$
corresponding to $M$.

\vskip0.3truecm\noindent  {\it Claim 3: When $M$ and $x,x'$ vary, $D_{M,x,x'}$ varies in a
$2$-dimensional system ${\cal D}$ of curves on $P$ with no base point. There is only one curve of ${\cal D}$ containing two general points of $P$,
 i.e. ${\cal D}$
has index $1$, hence it is an involution.} \vskip0.3truecm

Let $M$ be a general curve in ${\cal M}$, thus corresponding to a general point $m$ of $Pic^0(S)$. Of
course $M$ belongs to a $1$-dimensional system of curves $D_{M,x,x'}$, when $x,x'\in M$ are conjugated by
$\tau$. This proves that ${\cal D}$ is $2$-dimensional. A base point of ${\cal D}$ would correspond to
a curve $\overline M$ of ${\cal M}$ which belongs to $D_{M,x,x'}$ for the general curve $M$ and every pair of
points $x,x'$ conjugated in $\tau$ on $M$. But then $\overline M$ would have would have every pair of
points $x,x'$ on $M$ conjugated in $\tau$ in
common with $M$, a contradiction. The final assertion
follows by claim 1. \vskip0.3truecm

\vskip0.3truecm\noindent  {\it Claim 4: ${\cal D}$ is not a linear system.} \vskip0.3truecm

Suppose ${\cal D}$ is a linear system. Consider the morphism $\phi_{\cal D}: P\to \bbP^2$
determined by ${\cal D}$, which has degree at least $2$. This means that, given a general curve $M$,
corresponding to $m\in P$,  there is a curve $M'\not= M$ corresponding to $m'\in P$ with $m'\not= m$, such
that for every curve $D\in {\cal D}$ containing $m$, it also contains $m'$. Therefore for every pair of
points $x,x'$ conjugated in $\tau$ on $M$ the curve
$D_{M,x,x'}$, which contains $m$, also contains $m'$, and this implies that $M'$ has $x$ and $x'$ in common
with $M$. As $x,x'$ vary on $M$ staying conjugated in $\tau$, we see that $M$ and $M'$ have infinitely many
points in common, a contradiction.

\vskip0.3truecm\noindent  {\it Claim 5: ${\cal D}$ is not composed with a pencil.} \vskip0.3truecm

Suppose ${\cal D}$ is composed with a pencil. By the very definition of a family ${\cal M}_{M,x,x'}$, we
have that the general curve $D_{M,x,x'}$, if reducible, has all of its components containing the point
$m\in P$ corresponding to $M$. On the other hand, by the definiton of a system composed with a pencil,  the general curve
of such a system may have a singular point only at the base points of the pencil, which are fixed. Hence the general
curve of a system composed with a pencil is not singular at a moving point. Thus we see that $D_{M,x,x'}$ must be
irreducile. Since we are assuming that ${\cal D}$ is composed with a pencil, this would imply that ${\cal D}$ itself is a
pencil, which contradicts the fact that ${\cal D}$ has dimension $2$.\vskip0.3truecm

In conclusion claims 4 and 5 above contradict Castelnuovo-Humbert theorem above, which concludes our
proof in case $F=0$. \vskip0.3truecm

Next we consider the case $F\neq 0$. By lemma \ref{con}, $F$ is $1$-connected, $F\cdot M_{\eta}=2$ and  the linear system
$|2K_S|$ cuts out on $M$ a base point free linear series
$g^r_{2g}$, with
$r\geq g-1$. Suppose that this series is composed with an involution $\tau:=\tau_M$ of degree
$\delta\geq 2$. Then we must have $2g\geq \delta(g-1)$. This yields $\delta=2$.
Otherwise we would have $g\leq 3$, whereas $M^2\geq 3$, (see lemma \ref{paracanonicalcurve}), which implies $g\geq 5$.

If $r=g$, then
 $M$ is hyperelliptic and $|2K_S|$ cuts on $M$ the $g$-tuple
multiple of the $g^1_2$. In this situation, claim 2 above still holds. On the other hand, by
arguing as in claim 1 above, we see that, if $x, x'$ are two general points on $M$ conjugated in the
hyperelliptic involution, then $M$ is the unique curve in ${\cal M}$ containing them. Putting these
two informations together, we reach a contradiction.\par

If $r=g-1$ either $M$ is
hyperellitpic and we can argue as before, or $\phi(M)$ has arithmetic genus $1$, and then we can
argue as in the case $F=0$.\qed

\end{pf}\vskip0.5truecm

Now we are ready to give the:

\paragraph{Proof of proposition \ref{main5}} Proposition \ref{k8} is the statement for $K_S^2=8$ so
we can assume that
$K_S^2\leq 7$. Then, by lemma
\ref{paracanonicalcurve}, the general  curve $M:=M_\eta$ in ${\cal M}$ is irreducible and, by proposition \ref{grado},
$\phi$ is birational on $M$. Set $M'=M_{-\eta}$. Since $M'$ is also a general curve in ${\cal M}$, $\phi$ is also
birational on $M'$. Let $x\in M$ be a general point and let $x'\notin M$ be another point of $S$ such
that $\phi(x)=\phi(x')$. By the
generality of
$x\in M$, the point
$x'$ is also a sufficiently general point on $S$, hence it does not lie on $F$. Since $M+M'+2F \in |2K_S|$ then $x'\in M'$.
Again by the generality of $x'$ and of $M'$,  there is no other point $x"\in M'$ such that $\phi(x")=\phi(x)=\phi(x')$. So the
degree of $\phi$ has to be $2$.\qed

\bigskip

\section{The bicanonical involution}\label{sete}
 Let $S$ be a surface with $p_g=q=2$ presenting the non-standard case. By proposition \ref{main5} the bicanonical map
$\phi:S\to
\Sigma$    has degree $2$.
\par In general if the bicanonical map of a surface $S$ has degree $2$    we
can consider the bicanonical involution $ i: S\to S$.

The involution $i$ is biregular,
since $S$ is minimal of general type, and the fixed locus of $i$ is the
union of a smooth curve $R'$ and of isolated points $P_1,\ldots,
P_t$ . Let $\tilde{\Si}$ be the quotient of $S$ by $i$ and let
$p: S\to \tilde{\Si}$ be the projection onto the quotient.
The surface $\tilde{\Si}$ has nodes at the points
$Q_i:=p(P_i)$,
$i=1\ldots
t$, and is smooth elsewhere. Of course the bicanonical map of $S$ factors through $p$.
\par If $R'\neq \emptyset$, the image via
$p$ of $R'$ is a smooth curve $B''$ not containing the singular points $Q_i$, $i=1\ldots t$.

 Let now $f: V\to S$ be the
blow-up of $S$ at $P_1\ldots P_{t}$ and set $R=f^*R'$,
$E_i=f^{-1} (P_i)$, $i=1\ldots t$. The involution $i$ induces a
biregular involution $\tilde{i}$ on $V$ whose fixed locus is $R+\sum
E_i$. The quotient $W=V/<\tilde{i}>$ is smooth and one has a commutative
diagram:
 \begin{equation}
 \renewcommand{\arraystretch}{1.3}
 \begin{array}{ccc}
 \hphantom{\scriptstyle{g'}\ } V &\buildrel f\over \longrightarrow{} & S
 \hphantom{\ \scriptstyle{g}} \\
 {\scriptstyle{\pi}\ }\big\downarrow && \big\downarrow \ \scriptstyle{p} \\
 \hphantom{\scriptstyle{q}\ } W & \buildrel g\over \longrightarrow{} & \tilde{\Si}
 \hphantom{\ \scriptstyle{p}}
\end{array}
\end{equation}
where $\pi: V\to W$  is the  projection onto the quotient and
$g: W\to\tilde{\Si}$ is the minimal desingularization map. Of course
also the bicanonical map of $V$ factors through $\pi$. Notice  that
$A_i:=g^{-1} (Q_i)$ is an irreducible
$(-2)-$curve for
$i=1\ldots
t$. The map $\pi$ is flat, since it is finite and $W$ is smooth. Set $B'=g^*B''$.
 Thus there exists a line bundle $L$ on $W$ such that $2L\equiv B:= B'+\sum
A_i$ and $\pi_*\calO_V=\calO_W\oplus L^{-1}$.  ${\calO}_W$ is the invariant and $L^{-1}$ the
antiinvariant part of $\pi_* {\calO}_V$ under the action of $\tilde{i}$.  Since $\pi$ is a
double cover,
 the invariants of $V$ and $W$
relate by:
\begin{equation}\label{formule}
 \begin{aligned} K_V^2&=2(K_W+L)^2,\\
 \chi(\calO_V)&=2\chi({\calO}_W)+\frac{1}{2}L\cdot (K_W+L),\\
 p_g(V)&=p_g(W)+h^0(W,\calO_W(K_W+L)).
 \end{aligned}
 \end{equation}

Since $V$ is the blow-up of $S$ at $t$ points, $\chi(\calO_S)=\chi(\calO_V)$ and $K_S^2=K_V^2+t$.
In this case, because we are considering double covers through which the bicanonical map factors, we can be more precise:

\begin{prop} \label{KW} Let  $S$ be a minimal surface of general type with   $p_g(S)\geq 1$ and  bicanonical map of
degree $2$. Then, keeping the above notation, one has:

\vskip0.2truecm\noindent  i)  $h^0(W, {\calO}_W(2K_W+ L))=0$, $h^0(W, {\calO}_W(2K_W+ B))=h^0(S, {\calO}_S(2K_S))$;
\vskip0.2truecm\noindent  ii) either  $p_g(W)=0$ and $h^0(W, {\calO}_W(K_W+ L))=p_g(S)$\par  or
$\quad p_g(W)=p_g(S)$ and
$h^0(W, {\calO}_W(K_W+ L))=0$;
\vskip0.2truecm\noindent  iii) $|2K_V|=\pi^*|2K_W+B'|+\sum_i E_i$, $f^*|2K_S|=\pi^*|2K_W+B'|$ and
furthermore ${\calO}_W(2K_W+B')$ is nef and big;
\vskip0.2truecm\noindent  iv) $(2K_W+B')^2=2K_S^2$;
\vskip0.2truecm\noindent  v) $\chi({\calO}_W(2K_W+L))=0$;
\vskip0.2truecm\noindent  vi) $K_W\cdot (K_W+L)=\chi({\calO}_W)-\chi({\calO}_S)$.

\end{prop}
\begin{pf}
\vskip0.2truecm\noindent i), ii) By the projection formulas for double covers, one has
$$H^0(V,
{\calO}_V(K_V))=H^0(W,{\calO}_W (K_W))\oplus H^0(W, {\calO}_W (K_W+L))$$
 and
$$H^0(V,
{\calO}_V (2K_V))=H^0(W,{\calO}_W (2K_W+B))\oplus H^0(W, {\calO}_W (2K_W+L)).$$
\par

In both the above decompositions, the first summand is the invariant, the second the anti-invariant, part by the action of the
involution $\tilde i$. The fact that the bicanonical map of $V$ factors through $\pi$ implies the vanishing of one of the two
summands in each of the decompositions. Thus assertion ii) immediately follows. Since $p_g(S)\geq 1$, either the invariant or
the anti-invariant part of $H^0(V,
{\calO}_V(K_V))$ is non-zero. Hence the invariant part of $H^0(V,
{\calO}_V(2K_V))$ is certainly non-zero, and therefore i) also holds.

\vskip0.2truecm iii) Recall that $B=B'+\sum A_i$. Part i) implies that
$|2K_V|=\pi^*|2K_W+B|$. Since $|2K_S|$ is base point free (see \cite {ci}), then the fixed part of
$|2K_V|$ is $2\sum_i E_i$. More precisely, one has $|2K_V|=f^*|2K_S|+2\sum_i E_i$.  Thus one has $f^*|2K_S|=\pi^*|2K_W+B'|$
and therefore $ {\calO}_W(2K_W+B')$ is nef and big because ${\calO}_S(2K_S)$ is nef and big.

\vskip0.2truecm iv) follows immediately from $f^*|2K_S|=\pi^*|2K_W+B'|$ because $f^*(2K_S)^2=4K_S^2$ and $\pi$
is a double cover.

 \vskip0.2truecm v)  Since  $2(K_W+L)\equiv (2K_W+B')+\sum A_i$ and ${\calO}_W(2K_W+B')$
is nef and big by iii),  we can apply  the Kawamata-Viehweg's vanishing theorem to the divisor $K_W+L$ (see \cite{kawa},
corollary 5.12, c), pp. 48-49) obtaining:

 $$h^i(W,{\calO}_W (2K_W+L))=0\quad i=1,2.$$ By i) $h^0(W, {\calO}_W(2K_W+ L))=0$, thus
$\chi({\calO}_W(2K_W+L))=0$.

\vskip0.2truecm vi) By the Riemann-Roch theorem and by the formulas (\ref{formule}) we have:

 $$\chi({\calO}_W(2K_W+L))=\frac{1}
{2}(2K_W+L)\cdot (K_W+L)+\chi({\calO}_W)=$$
 $$=K_W\cdot (K_W+L)+ \frac{1}
{2}
L\cdot (K_W+L)+\chi({\calO}_W)=$$
$$=K_W\cdot (K_W+L)+\chi({\calO}_S)-\chi({\calO}_W).$$\medskip
Then the assertion follows from part v).
\qed  \end{pf}

\vskip0.5truecm If $S$ is a minimal surface of general type with $p_g=q=2$ and
bicanonical map of degree $2$, we can be more specific.
\begin{lem} \label{invariants} Let $S$ be a minimal surface of general type with $p_g=q=2$ for which the Albanese map is
surjective.  Suppose  the bicanonical map of $S$ has degree  $2$ and let
 $W$ be as above. Then either: \par

\noindent i) $p_g(W)=2,q(W)=2$,

or \par\noindent ii) $p_g(W)=0, q(W)=1$, \par
 or
\par\noindent iii)  $p_g(W)=2, q(W)=0$.
\end{lem}

\begin{pf} By ii) of proposition \ref{KW} we know that either $p_g(W)=2$ or $p_g(W)=0$. By the projection formulas for
double covers, one has
$$2=q(S)=H^1(V,{\calO}_V (K_V))=H^1(W,{\calO}_W (K_W))\oplus H^1(W,{\calO}_W( K_W+L))$$ and therefore $q(W)\leq 2$
 with equality holding if and only if $h ^1(W, {\calO}_W(K_W+L))=0$.\par Assume that $q(W)=2$. Then
$H^0(V,\Omega^1_V)$  is generated by two  1-forms $\omega, \omega'$ which are invariant under the bicanonical
involution and therefore $\omega\wedge\omega'$ is an invariant element of $H^0(V, \Omega^2_V)$. Since, by
corollary \ref{defranchis},
$\omega\wedge\omega'\not\equiv 0$,
$p_g(W)\neq 0$ and so $p_g(W)=2$.\par Assume now that $q(W)=1$. Then $H^0(V,\Omega^1_V)$ has invariant and
antiinvariant subspaces both of dimension $1$. If $\omega^+$ and $\omega^-$ are generators of such subspaces, they
form a basis of $H^0(V,\Omega^1_V)$. Since, as before, $\omega\wedge\omega'\not\equiv 0$,  $\omega^+\wedge \omega^-$
is a nonzero antiinvariant element of $H^0(V,\Omega^2_V)$.  So $p_g(W)$ is not $2$ and therefore $p_g(W)=0$.
\par Suppose now that $q(W)=0$. Then $H^0(V,\Omega^1_V)$  is generated by two  1-forms $\omega, \omega'$ which are
antiinvariant under the bicanonical involution and therefore $\omega\wedge\omega'$ is an invariant element of
$H^0(V, \Omega^2_V)$. As in the preceding paragraphs we conclude that $p_g(W)=2$.
 \qed  \end{pf}

We keep the same assumptions as in lemma \ref{invariants}, and we analyse the possibilities given by the lemma.

\begin{lem} \label{q2}
The case $q(W)=2$ cannot occur.
\end{lem}
\begin{pf} Suppose otherwise. By proposition \ref{KW}, vi) we have $K_W\cdot (K_W+L)=0$ and so $K_W\cdot
(2K_W+B')=0$.  Therefore, since $|K_W|$ is  a pencil
we get a contradiction to the fact that $2K_W+B'$ is nef and big, (see proposition  \ref{KW}, iv)).
\qed  \end{pf}

\begin{lem} \label{q1}
 Keep the assumptions in lemma \ref{invariants} and assume furthermore  that $S$ has no genus $2$ pencils. Then the
case
$q(W)=1$ does not occur.
\end{lem}

\begin{pf}  We notice first that $k(W)<0$ and thus $W$ is a ruled surface. In fact suppose otherwise. Then  some
multiple of
$K_W$ is an effective divisor. By proposition
\ref{KW}, vi) we have
$K_W\cdot (K_W+L)=-1$, and so $K_W\cdot (2K_W+B')<0$, which contradicts $2K_W+B'$ being nef and big.
\par

In this case  we have, by proposition \ref{KW}, ii), $h^0(W,{\calO}_W(K_W+L))=2$ and thus we can write
$|K_W+L|=|Y|+Z$, where $|Y|$ is the moving part and  $Z$ is the fixed part. Since for each
$(-2)-$curve $A_i$ we have $A_i\cdot (K_W+L)=-1$, then $Z\neq 0$.  Notice that  $\pi^*(|Y|)$ is exactly the moving part of
$|K_V|$ and therefore by proposition \ref{irreducible} the general curve $Y$ in
$|Y|$ is irreducible. Furthermore, since $W$ is not rational, and $|Y|$ is a linear system of dimension
$1$, the geometric genus of a general curve $Y\in |Y|$ is at least $1$.

\vskip0.3truecm\noindent {\it Claim 1 :
\noindent  for every effective, non-zero divisor $N< K_W+L$, one has
$h^0(N,{\calO}_N)+p_a(N)\leq 2$.
}\vskip0.3truecm
By the Riemann-Roch theorem, we have $$h^0(W,{\calO}_W(K_W+N))+ h^2(W,{\calO}_W(K_W+N))=$$ $$=\frac{1}
{2}(K_W\cdot
N+N^2)+h^1(W,{\calO}_W(K_W+N)). \quad\quad (*)$$
Now notice  that $h^0(W,{\calO}_W(K_W+N))=0$. If not, since
$N< K_W+L$, we would have
$h^0(W,{\calO}_W(2K_W+ L))\neq 0$, a contradiction to proposition \ref {KW}.  Since $N$ is effective
we have also $h^2(W,{\calO}_W(K_W+N))=0$,  and so $(*)$ can be written as
$$p_a(N)-1+h^1(W,{\calO}_W(K_W+N))=0.$$

 Since, by  duality,  $h^1(W,{\calO}_W(K_W+N))=h^1(W,{\calO}_W(-N))$ and  one has $h^1(W,{\calO}_W(-N))\geq
h^0(N,{\calO}_N)-1$ for any effective divisor $N$,  we obtain $p_a(N)-1 +
h^0(N,{\calO}_N)-1\leq 0$, proving the claim.\par

\vskip0.3truecm\noindent {\it Claim 2 : If  $T$ is a general ruling of $W$, then $Y\cdot T=1$. }\vskip0.3truecm

 As we  noticed already,  the geometric genus of a general curve $Y\in |Y|$ is at least
$1$, and of course $h^0(Y,{\calO}_Y)=1$. By claim 1 we conclude that
 $Y$ is smooth and elliptic.  Claim $1$ implies also that each irreducible component
$\theta$ of
$Z$ is rational and  such that $\theta\cdot Y\leq 1$.  Consider the pencil $|Y|$.  By the Riemann Roch
theorem $$h^0(W, \calO_W(Y))=Y^2+h^1(W,
\calO_W(Y))$$ and so $0\leq Y^2\leq 2$.  We claim that
$|Y|$ has no multiple fibres. If
$Y^2>0$ the claim is trivial,  since $Y^2\leq 2$. Assume $Y^2=0$ and notice that $Y\cdot L>0$ because otherwise we would have a
pencil of
 curves of genus $1$ on $V$, which is impossible. Hence $Y\cdot Z=Y^2+Y\cdot Z=
Y\cdot (K_W+L)=Y\cdot L>0$ and thus there exists  an irreducible curve $\theta$ in $Z$ such that $Y\cdot
\theta=1$. So  also in this case
$|Y|$ has no multiple fibres. \par

We can  consider now the relatively minimal fibration $h:\tilde W\to {\bbP^1}$ associated to $|Y|$, i.e. we blow up the base
points of $|Y|$, if any, and contract the $(-1)$-curves contained in fibres of
$|Y|$. Since $|Y|$ has no multiple
fibres and $\chi({\calO}_{\tilde W})=1$, we have by \cite{bpv}, corollary  V.12.3, pg. 162, $K_{\tilde W}\equiv
-2F$, where $F$ is a general fibre of
$h$.
\par Let now $T$ be a general ruling of $W$ and $\tilde T$ the corresponding ruling of $\tilde W$. Since $K_{\tilde
W}\cdot {\tilde T}=-2$, we conclude that $F\cdot {\tilde T}=1$ and therefore also $Y\cdot T=1$ proving claim $2$.

\vskip0.3truecm

Now we can finish our proof. Let $T$ be the general ruling of $W$. Since each component of $Z$ is rational, $TZ=0$, and
so  we have
$(K_W+L)\cdot T=(Z+Y)\cdot T= Z\cdot T+Y\cdot T=1$. Since $K_W\cdot T=-2$, we have $L\cdot T=3$. This implies that, by
pulling back to $V$ the ruling of $W$, we obtain a  pencil of curves of genus $2$, against our hypothesis.
 \qed  \end{pf}

Finally we come to the  case $q(W)=0$.  \par

\begin{prop}\label {nodes}  Keep the assumptions as in lemma \ref{invariants} and assume furthermore  that $S$ has no genus
$2$ pencils. If $q(W)=0$ then $B'=0$, $W$ is a minimal surface of general type with $p_g(W)=2$,
$K_W^2=2$ and $p: S\to \tilde\Sigma$ is ramified only at $20$ nodes of $\tilde\Sigma$. Furthermore, if $C$ and $C'$ are
the general  curves in
$|K_W|$ and
$|K_S|$ respectively, $C$ and $C'$ are smooth irreducible non hyperelliptic.

 \end{prop}

\begin{pf} We keep the notation as in the beginning of the section. \par
Let $a:S\to A:=Alb(S)$ be the Albanese map.
We can define a morphism $\tilde a:\tilde\Sigma\to A$ by associating to each point $x\in \tilde\Sigma$ the sum
of the Albanese images of the two points in the cycle $p^*(x)$. Since $q(\tilde\Sigma)=0$ this map is constant and, up to a
translation, we may assume that its image is the point $0\in A$. Hence if $p^*(x)=y_1+y_2$ we have $a(y_1)=-a(y_2)$. Thus we can
define a morphism $\alpha:\tilde\Sigma\to K(A)$, where $K(A)$ is the Kummer surface of $A$, by associating to
$x\in\tilde\Sigma$ the point $y\in K(A)$ corresponding to $a(y_1)=-a(y_2)$.

Given any point $x_0$ in the branch locus, we have
$p^*(x_0)=y_0+y_0$, so $a(y_0)$ is a 2-torsion point in $A$.  In particular the ramification divisor
$R$ must be contracted by the Albanese map and so also $B''$ is contracted by $\alpha$.
Notice that $K_{\tilde\Sigma}=\alpha^*(K_{K(A)})+D$, where $D$ is the divisor where the differential of $\alpha$ drops rank,
in particular $D$ contains all the curves contracted by $\alpha$. Since $K(A)$
is a $K3$ surface, we see that there is an effective canonical divisor on
${\tilde\Sigma}$ containing the smooth curve $B''$. Hence also $K_W$
can be written as
$B'+\Delta$, where
$\Delta$ is an effective divisor.\par

  Notice that by the classification of surfaces  $W$ is either
elliptic or of general type. \par
 Let $|K_W|=|Y|+Z$, where $Z$ is the fixed part and $|Y|$ the movable part of $|K_W|$. Since
$p_g(W)=2$, the system $|Y|$ is a pencil. Since $W$ is regular, by Bertini's theorem the general curve
of $|Y|$ is irreducible.

 Remember that  the bicanonical map of $V$ has degree $2$ to its
image, and factors through $\pi$ and through the map defined by the linear system $|2K_W+B'|$ on $W$. This
implies that the linear series cut out by $|2K_W+B'|$ on the general curve $Y\in |Y|$ determines a birational
map on $Y$. In particular it has projective dimension at least $2$. Since $W$ is not ruled, one has $(2K_W+B')\cdot Y\geq
3$, with equality being possible only if $g(Y)=1$, which in turn is only possible if $Y^2=K_W\cdot Y=0$.

By the formulas (\ref{formule}) and by proposition
\ref{KW},  we have
 $K_W\cdot (K_W+L)=2$, hence $(2K_W+B')\cdot K_W=4$.

Since $2K_W+B'$
is nef we have
$(2K_W+B')\cdot Y\leq (2K_W+B')\cdot K_W=4$. Since $Y$ is nef, one has  $B'\cdot Y\geq 0$, hence we obtain
$K_W\cdot Y\leq 2$. By the adjunction formula $Y\cdot Z$ is even, hence
either $Y\cdot Z=0$ or $Y\cdot Z=2$. \par
On the other hand we have seen above that we can write $K_W= B'+\Delta$, where $\Delta$ is an effective divisor, and
so
$2K_W+B'= 3B'+2\Delta$, hence
$3\leq 3B'\cdot Y+2Y\cdot \Delta\leq 4$, so either
$B'\cdot Y=0$ or $B'\cdot Y=1$ and $g(Y)=1$. This is impossible because then $Y^2=Y\cdot K_W=0$, thus $0=Y\cdot K_W=Y\cdot
B'+Y\cdot \Delta= 1+Y\cdot \Delta$ and the nef divisor
$Y$ would be such that $Y\cdot \Delta=-1$, a contradiction. Thus the only possibility is $Y\cdot B'=0$,
$Y\cdot \Delta=2$ and therefore
$K_W\cdot Y=2$.
Since $Y\cdot Z$ is even and non negative, we have that either $Y^2=0, Y\cdot Z=2$ or $Y^2=2, Y\cdot Z=0$. \par

In the
first case $Y^2=0, Y\cdot Z=2$, we get
 ${\calO}_Y(2K_W+B')\simeq {\calO}_Y(2K_Y)$. This  is impossible because in this case $|Y|$ is a genus $2$ pencil and
so
$|2K_W+B'|$ would determine a non-birational map on
$W$.

If $Y^2=2, Y\cdot Z=0$, then we have $2K_W+B'\equiv 2Y+(2Z+B')$ and $2Y\cdot (2Z+B')=0$. Since $2K_W+B'$ is nef and big, the
only possibility is that $2Z+B'= 0$. So $B'=Z=0$, hence $K_W\equiv Y$ is nef and therefore $W$ is minimal. Moreover
$K_W^2=Y^2=2$. Furthermore $2K_W+B'=2K_W$ and, by proposition \ref{KW}, iv) we have $K_S^2=4$.  In addition, by
the formulas (\ref{formule}) and by proposition \ref{KW},  we have
  $(K_W+L)^2=-8$ and so $K_V^2=-16$. Hence $t=16+K_S^2=20$, where $t$ is, as before, the number of isolated fixed points of
the bicanonical involution. Thus
$p$ is ramified exactly over 20 nodes.\par

By the above the general curve $C$ in the linear system  $|K_W|=|Y|$ is irreducible and non hyperelliptic, because the
bicanonical map of $W$ is birational. Since $K_W^2=2$, and $|K_W|$ is a rational pencil, $C$ is necessarily smooth.
The assertion for the general curve $C'$ in $|K_S|$ is then obvious. \qed
\end{pf}

\section{ The  main theorem.}\label{main}
In the previous sections we saw that if the bicanonical map $\phi$ of a surface $S$ with $p_g=q=2$ is not birational and $S$
has no pencil of curves of genus $2$ then $\phi$ has
degree $\sigma=2$ and we have described in proposition \ref{nodes} some properties of the quotient of $S$  by the involution
induced by the bicanonical map.\par

In this section we will classify these quotients. Let us start by presenting an example, which was first
pointed out by F. Catanese  (cf.  \cite{ci},  example (c), page 70, and remark 3.15, page 72).

\begin{ex} \label{Ciro}
 {\rm  Let $A$ be an abelian
surface with an irreducible symmetric principal polarization $\Theta$, and suppose that $A$ contains no elliptic curves. Let
$h\colon S\to A$
be the double  cover branched on a smooth divisor $B\in |2\Theta|$ so
that $h_*\calO_S=\calO_A\oplus \calO_A(-\Theta)$. Since
$K_S=h^*(\Theta)$, the invariants
of the smooth surface $S$ are $p_g(S)=2$, $q(S)=2$, $K^2_S=4$. Notice that the map $h\colon
S\to A$ factors through the Albanese map $a: S\to Alb(S)$. Since $h$ has degree $2$ and $Alb(S)$ is a surface,
we see that that $Alb(S)\simeq A$. In addition we observe that $S$ has no genus $b$ pencil of curves of genus $2$. Indeed, by
lemma \ref {fibration}, ii) and by the assumption that $A\simeq Pic^0(S)^*$ contains no elliptic curve, one should have $b=2$,
and by part i) of the same lemma we would find $K_S^2\geq 8$, a contradiction.\par

Remark now that $B$ is
symmetric with respect to the involution $j$ of $A$ determined by the multiplication by $-1$.
Hence $j$ can be
lifted to an involution $i$ on $S$ that
acts as the identity on $H^0(S, \calO_S(K_S))$. We denote by $p\colon S\to \tilde \Sigma :=S/<i>$
the
projection onto the quotient.  We observe that
$p_g(\tilde \Sigma)=2$, $q(\tilde \Sigma)=0$, $K^2_{\tilde \Sigma}=2$ and the only singularities  of the surface
$\tilde \Sigma$ are $20$
nodes. Since $h^0(\tilde \Sigma,\calO_{\tilde \Sigma}(2K_{\tilde \Sigma}))=\chi(\calO_{\tilde
\Sigma})+K^2_{\tilde
\Sigma}=5=h^0(S,\calO_S( 2K_S))$, the bicanonical map of
$S$ factors through $p\colon S\to {\tilde\Sigma}$. Since $S$ has no pencil of curves of genus $2$, we have the
situation described in proposition \ref{nodes}.

\par For the sake of completness, we want to point out the following  alternative description of
$\tilde\Sigma$. One embeds, as usual, the Kummer surface $Kum(A)$ of $A$ as a
quartic
surface in $\pp^3=\pp(H^0(A,2\Theta)^*)$. The surface $\tilde\Sigma$ is a double  cover
of
$Kum(A)$  branched along the smooth plane section $H$ of $Kum(A)$
corresponding
to $B$ and on $6$ nodes, corresponding to the six points of order 2 of $A$
lying on $\Theta$. The ramification divisor $R$ of ${\tilde \Sigma}\to Kum(A)$ is a canonical
curve isomorphic to $H$, and thus it is not hyperelliptic.} \end{ex}

\begin{rem}\label{deg4}  {\rm The same construction can also be done with a reducible  polarization $\Theta$  on
$A$. Then
$A$  is  isomorphic to the product
$E_1\times E_2$ of two elliptic curves and the surface $S$ constructed as above  has two elliptic pencils of
curves of genus
$2$ curves. In this case the bicanonical map of $S$ has  degree $4$ ( see \cite{xiao}, Thm. 5.6).} \end{rem}
\vskip0.5truecm

We are finally going to prove our classification theorem:
\begin{thm} \label{maint} Let $S$ be a minimal surface of general type with
$p_g=q=2$, presenting the non standard case. Then $S$ is as in
example \ref{Ciro}.\end{thm}

For the proof we need a
preliminary lemma and  some notation. Let $X$ and $Y$ be smooth, projective
surfaces and
$f: X\to Y$ be a surjective map. Let $R$ be the ramification curve on $X$, i.e. the subscheme of $X$ where $f$
drops rank. Let
$C$ be a smooth, irreducible curve on
$X$ non contained in $R$. Set $\Gamma:=f(C)$ and $f^*(\Gamma)=C+D$. Notice that $C$ and $D$ have no common
component. For every point $p\in C$, denote by $r_p$ [resp. by $d_p$] the coefficient of $p$ in the divisor
cut out on $C$ by $R$ [resp. by $D$]. Set $\delta_p=r_p-d_p$ and $p':=f(p)$. Then:

\begin{lem} \label{delta} With the above notation, if $\Gamma$ is smooth at
$p'$, then $\delta_p\geq 0$. \end{lem}

\begin{pf} Put local coordinates $(s,t)$ centered at $p'$ in such a way that
$\Gamma$ has equation $t=0$. Put local coordinates $(x,y)$ centered at $p$ in such a way that $C$ has
equation $x=0$ and $\phi(x,y)=0$ is the equation of $D$. Then $f$ has local equations $s=\psi(x,y)$
and $t=x\phi(x,y)$. Therefore $R$ has equation:

$$\phi{\frac{\partial \psi}{\partial y}}+x {\frac{\partial (\phi,\psi)}{\partial(x, y)}}=0$$

\noindent whence the assertion immediately follows. \qed
\end{pf}

Now we can prove our classification theorem:

\paragraph{Proof of theorem \ref{maint}} The main step in our proof is to show that the Albanese map $a: S\to
A:=Alb(S)$ has degree $\nu=2$. This is what we are going to prove first.\par

As we saw in section \ref{sete}, the bicanonical map of $S$ factors through the degree 2
finite cover
$p:S\to \tilde\Sigma$ branched only at the 20 nodes of $\tilde\Sigma$. By proposition \ref{KW}, iii),
$ K_{\tilde\Sigma}$ is a nef and big line bundle on $\tilde\Sigma$. More precisely, from proposition \ref{nodes}
it follows that $|K_{\tilde\Sigma}|$ is a pencil with no fixed component and with 2 base points which do not
occur at any of the nodes of $\tilde\Sigma$. Hence $(p:S\to
\tilde\Sigma, K_{\tilde\Sigma})$ is a good generating pair in the sense of \cite{cpt}.\par

Let $C$ be a general curve in
$|K_{\tilde\Sigma}|$ and let $C':=p^*(C)$. Since $C$ does not contain any of the nodes of $\tilde\Sigma$, then
$p:C'\to C$ is an \'etale double cover. Theorem (6.1) of [CPT] yields that the Prym variety $P:=Prym(C',C)$
related to the double cover $p:C'\to C$ is isomorphic to the Albanese surface $A$. Therefore $A$ is
principally polarized, and we denote by $\Theta$ its principal polarization. Furthermore, after
having identified $A$ with $P$, the Abel-Prym map $\alpha: C'\to P$ coincides, up to
translation, with the restriction to $C'$ of the Albanese map $a: S\to A$. Notice that $C'$ is not
hyperelliptic and set $\Gamma:=a(C')$. By the results in \cite{lb}, chapter 12, the map $a_{|C«}: C'\to
\Gamma$ is an isomorphism and therefore $\Gamma$ is smooth. Furthermore $\Gamma$ is in the class of $2\Theta$ by Welters'
criterion (see again \cite{lb}, chapter 12).\par

 Let us set $a^*(\Gamma)=C'+D$ and let us denote by $R$ the ramification curve of
$a$. By lemma \ref{delta}  we have
$K_S\cdot D=C'\cdot D\leq C'\cdot R=C'\cdot K_S=4$, with equality holding if and only if $D\sim K_S$. By the
index theorem we have
$D^2\leq 4$. Thus
$\nu\cdot 8=\nu\cdot (2\Theta)^2=a^*(\Gamma)^2=(C'+D)^2\leq 16$. This proves that $\nu=2$ and, in addition,
that $D\sim K_S$.\par

Now we can finish our proof by showing that the branch curve $B$ of $a: S\to A$ is a divisor in the class of
$2\Theta$. This immediately follows from the fact that $16=2R\cdot (C'+D)=2B\cdot \Gamma$, hence $B\cdot
\Theta=4$, so that $B$ is numerically equivalent to $2\Theta$. \qed

\bigskip

\bigskip

\begin{tabbing}  Universit\`a di Roma Tor Vergata xxxxxxxxxxxxxxxx\=Universidade de Lisboa  \kill Ciro
Ciliberto \> Margarida Mendes Lopes\\
 Dipartimento di Matematica \> CMAF\\
 Universit\`a di Roma Tor Vergata \>Universidade de Lisboa \\ Via della Ricerca Scientifica \>
Av. Prof. Gama Pinto, 2 \\ 00133 Roma, ITALY \> 1649-003 Lisboa, PORTUGAL\\
cilibert@mat.uniroma2.it \> mmlopes@lmc.fc.ul.pt
\end{tabbing}

\end{document}